\documentclass[12pt]{article}
\usepackage{latexsym,amssymb,amsmath,url}

\def\N{\mathbb{N}}
\def\Q{\mathbb{Q}}
\def\Z{\mathbb{Z}}
\def\R{\mathbb{R}}
\def\C{\mathbb{C}}

\def\eproof{\hfill{$\Box$}\bigskip}

\def\ds{\dots}
\def\sus{\subset}
\def\al{\alpha}
\def\be{\beta}

\def\cc{\colon}
\def\ep{\varepsilon}

\newtheorem{thm}{Theorem}[section]
\newtheorem{prop}[thm]{Proposition}

\newtheorem{defi}[thm]{Definition}

\begin{document}
\title{Effective formulas for linear recurrence sequences of integers\footnote{This text did not
make it to a conference, one of the reasons being ``Finally, while I personally find the results most interesting, 
I'm not confident that there would be an appeal to a broad audience at $\ds$''. For the 
interest of possible readers I put it in its rough form (sorry for that) at least to arXiv. I hope to produce more polished and 
complete version soon.}}
\author{Martin Klazar\footnote{Department of Applied Mathematics, Faculty of Mathematics and Physics, Charles University, 
Malostransk\'e n\'am. 25, 118 00 Praha 1, Czechia, {\tt klazar@kam.mff.cuni.cz}}}

\maketitle

\begin{abstract}
We propose a new definition of effective formulas for problems in enumerative combinatorics. 
We outline the proof of the fact that every linear recurrence sequence of integers has 
such a formula. It follows from a lower bound that can be deduced from the Skolem--Mahler--Lech theorem and the Subspace 
Theorem. We will give details of this deduction that is due to P.~Corvaja in the full version of this extended abstract.
\end{abstract}

\section{Introduction}

This is an extended abstract of a roughly twice as long article that in turn is a revised digest from the longer work \cite{klaz_I} of the author;
we therefore apologize for possible incoherence. Our original contribution is twofold. We propose a novel definition of
effective formulas for enumerative problems, and apply it to computation of terms in {\em linear recurrence sequences} in $\Z$. 
This may seem as resolved long time ago, but we explain that an unconditional algorithm computing terms in 
a general integral linear recurrence sequence that is effective according to our definition, became possible only after about 1981. 
In the first section we give definitions, state our results, and outline the algorithm. 
The second section contains the algorithm with outlined proofs. In the last third section we mention the 
Skolem--Mahler--Lech theorem and the Subspace Theorem that are needed for the proof of a lower bound, crucial for the algorithm.

An {\em enumerative problem} is a computable function
$$
f\cc\N\to\Z\;.
$$
Here $\N=\{1,2,\ds\}$ are the natural numbers and $\Z$ denotes the ring of integers. We postulate the codomain $\Z$, 
and not the more common $\N$ or $\N_0=\{0,1,\ds\}$, in order that a general integral linear recurrence sequence be an 
enumerative problem. A {\em linear recurrence sequence} is a function 
$f\cc\N\to\Z$ such that for some $k\in\N_0$ numbers $a_0,a_1,\ds,a_{k-1}\in\Z$ with $a_0\ne0$ and every $n\in\N$,
$$
f(n+k)=\sum_{i=0}^{k-1}a_if(n+i)\;.
$$
We call the number $k\in\N_0$ the {\em order of the recurrence}.
For $k=0$ we define the sum as $0$ and obtain the {\em zero sequence} with $f(n)=0$ for every $n\in\N$. Due to $a_0\ne0$ one can run the 
recurrence backwards for $n\le0$ and extend $f$ to $f\cc\Z\to\Q$. The {\em characteristic polynomial of the recurrence} is
$$
q_{\underline{a}}(x)=x^k-a_{k-1}x^{k-1}-\ds-a_1x-a_0,\ \underline{a}=(a_0,\,a_1,\,\ds,\,a_{k-1})\;.
$$
One of the most popular linear recurrence sequences is the {\em Fibonacci numbers}, given by $k=2$ and $f(1)=f(2)=a_0=a_1=1$. The characteristic polynomial is then $x^2-x-1$.

We define more generally a {\em linear recurrence sequence in $R$}, where $R$ is any integral domain, the {\em order of the recurrence}, and 
the {\em characteristic polynomial of the recurrence} in the same way by replacing $\Z$ with $R$. We need it for $R=\Q$ and $R=\overline{\Q}$ 
(the algebraic numbers).

Our new definition of effective formulas for enumerative problems is as follows. 

\begin{defi}
\label{our_defi}
A PIO formula for an enumerative problem $f\cc\N\to\Z$ is an algorithm that for two constants $c,d\in\N$ computes for each 
input $n\in\N$ the output $f(n)\in\Z$ in time at most
$$
c\cdot\big(\log(1+n)+\log(2+|f(n)|)\big)^d\;.
$$
\end{defi}

\noindent
Time here means the number of steps that a multitape Turing machine implementing the algorithm does when it computes $f(n)$ from the input $n$, 
where $n$ and $f(n)$ are represented as binary strings. Both logarithms are natural, and the shifts $1+\ds$ and $2+\ds$ 
remove inconvenient arguments $0$ and $1$. The definition could be simplified by using just one constant $d$, but this would also make it more rigid.

Definitions of effective formulas for problems in enumeration were proposed by H.\,S.~Wilf \cite{wilf} in 1982 and recently by I.~Pak \cite{pak}. 
Several other articles (usually building on \cite{wilf}) mention these definitions too, see \cite{klaz_I}, but for them it is not the main subject. 
We will not discuss two Wilf's definitions here and refer 
the reader directly to \cite{wilf} (some sources quote them imprecisely). Definition~\ref{our_defi} differs from those in \cite{wilf} and \cite{pak}
in two respects. It applies to any computable function from $\N$ to $\Z$, be it $\lfloor\log(1+\log n)\rfloor$ or the Ackermann function, while
\cite{wilf} restricts to functions in a growth range. More importantly, it takes into account the output complexity of the problem given 
by $f(n)$ while \cite{wilf} and \cite{pak} only work with the input complexity given by $n$.  For more discussion, examples and motivation 
related to Definition~\ref{our_defi} see \cite{klaz_I}.

The acronym PIO refers to the complexity class {\em polynomial input output}, see \cite{zoo}. By \cite{zoo} it was introduced by 
M.~Yanakakis \cite{yana}. Y.~Gurevich and S.~Shelah  introduced it (more precisely, a very similar complexity class) in \cite{gure_shel} 
independently later. The definition of the class in \cite{zoo} is similar to Definition~\ref{our_defi}, it consists of the maps 
$f\cc\{0,1\}^n\to\{0,1\}^m$ computable in time polynomial in $n$ and $m$.
The definition of the class in \cite{gure_shel}, where the acronym PIO means {\em computable in time polynomial in input or output}, 
concerns partial maps from $\Sigma_1^*$ to $\Sigma_2^*$ for two alphabets $\Sigma_i$, and is very similar to Definition~\ref{our_defi}.
The contribution of Definition~\ref{our_defi} is in relating the class PIO to enumerative problems, like linear 
recurrence sequences. We are not aware of any mention  prior to \cite{klaz_I} of this class in connection with problems 
in enumerative combinatorics. PIO in fact occurs rather rarely even in the computational complexity literature proper (see \cite{klaz_I} 
for a few more references). 
To finish the review of history of Definition~\ref{our_defi}, we mention that the present author stated it in a form in \cite{klaz} in 2010 
(unaware then of the class PIO) and that it was also briefly mentioned as a ``very good formula'' by J.~Shallit \cite[slide 3]{shal} 
in a lecture in 2016.

Definition~\ref{our_defi} is based on bit complexity. Other complexity measure is the algebraic complexity, dealing with numbers of arithmetic 
or other operations needed to compute functions. An algebraic complexity counterpart to Definition~\ref{our_defi} could be of an interest and 
we hope to present one in \cite{klaz_stir_conu}. 

The next section outlines proof of the following theorem, in fact of a stronger form given in Theorem~\ref{main_gene}.

\begin{thm}
\label{main}
Every linear recurrence sequence $f\cc\N\to\Z$ has a PIO formula. 
\end{thm}

\noindent
Why is this result not immediate? It is immediate for the Fibonacci numbers $F_n$ and for similar exponentially growing sequences. 
It is easy to prove by induction that $2^{n/2}\ll F_n\ll 2^n$ (on $\N$). Recall that for $f,g\cc M\to\C$, $M\sus\R$, the notation $f\ll g$ 
(on $M$) means that for some $c\in\N$ one has $|f(x)|\le c|g(x)|$ for every $x\in M$. The notation $f=O(g)$ (on $M$) has identical meaning. 
We introduce the notation $f=\mathrm{poly}(g)$ (on $M$) meaning that $g^c\ll f\ll g^d$ (on $M$)
for real constants $d\ge c>0$. If $M=\N$, we omit the qualifier ``(on $\N$)''. Thus for $f(n)=F_n$ the displayed expression 
in Definition~\ref{our_defi} is $\mathrm{poly}(n)$. There is the {\em simple algorithm} that directly applies the defining recurrence 
and computes $F_n$ also in time $\mathrm{poly}(n)$. It is therefore 
a PIO formula for the Fibonacci numbers. More generally, this algorithm is a PIO formula for any linear recurrence sequence $f\cc\N\to\Z$ with
an exponential lower bound $f(n)\gg c^n$ for a real $c>1$. By an easy extension of the inductive argument for $F_n$, every linear recurrence 
sequence has an exponential upper bound, but existence of an exponential lower bound is a problem. Many linear recurrence sequences lack it, 
for example
$$
f_1(n)=3^n+(-3)^n+n^4-2n=3^n+(-3)^n+(n^4-2n)\cdot 1^n
$$
(which is a linear recurrence sequence, see Section 2). The simple algorithm still computes  $f_1(n)$, of course, and may be viewed as effective 
by some criteria, but not by Definition~\ref{our_defi}. Indeed, for odd $n$ we have $f_1(n)=n^4-2n$ and the expression in Definition~\ref{our_defi} is $\mathrm{poly}(\log(1+n))$, but the simple algorithm takes much longer $\mathrm{poly}(n)$ time. Therefore it is not a PIO formula for $f_1(n)$. 

But in this particular example it is easy to fix. For even $n$ we still have an exponential lower bound, because $f_1(n)=2\cdot 3^n+n^4-2n$, 
and can compute $f_1(n)$ by the simple algorithm. For odd $n$ we compute $f_1(n)=n^4-2n$ faster in time $\mathrm{poly}(\log(1+n))$ 
by a direct evaluation of this integral polynomial, in three multiplications and one subtraction of two $O(\log(1+ n))$-bit integers. 
Let us call this better algorithm a {\em fast algorithm}. It splits the domain $\N$ of $f_1$ in two arithmetical progressions, 
of even numbers and of odd numbers, and computes $f_1(n)$ separately on each progression by a different (indicated) method.
It is clear that the fast algorithm is a PIO formula for $f_1(n)$. 

Section 2 treats in this way any linear recurrence sequence $f$. It turns out that for each $f$ there is a modulus $m\in\N$ 
such that on each residue class modulo $m$ the sequence $f$ either restricts to a rational polynomial with effectively bounded degree or 
has an exponential lower bound. In the former case we compute $f(n)$ by a direct polynomial evaluation, and in the latter case we use the simple algorithm. 
In the next section we give more details. The complexity bound of the algorithm rests on a deep result from the theory 
of linear recurrence sequences, an exponential lower bound on growth of non-degenerate sequences. We state it precisely in two equivalent forms.

Let $\overline{\Q}\sus\C$ be the field of algebraic numbers, consisting of the roots of rational polynomials, and 
$f\cc\N\to\overline{\Q}$ be a linear recurrence sequence in $\overline{\Q}$ with recurrence order $k$. 
If $k$ is minimum for the given sequence $f$, the recurrence is unique and we call the characteristic polynomial 
$$
q_f(x)=q_{\underline{a}}(x)=x^k-\sum_{i=0}^{k-1}a_ix^i
$$
the {\em characteristic polynomial of $f$}. Its roots are the {\em roots of $f$}. If the ratio of some two distinct roots 
of $f$ is a root of $1$, we say that $f$ is 
{\em degenerate}, and else that $f$ is {\em non-degenerate}. For example, if $f$ is the zero sequence then $q_f(x)=1$, $f$ 
has no root, and is non-degenerate. The above 
sequence $f_1$ has roots $-3,3$, and $1$ and is degenerate. 
Its restrictions to even, resp. odd, numbers (see Proposition~\ref{mjsec}) are non-degenerate
and have roots $3$ and $1$, resp. $1$. The next lower bound, partially proven by A.~van der Poorten in 1981, 
is crucial for the complexity analysis of the algorithm in Section 2.

\begin{thm}
\label{lower_b}
Suppose that $f\cc\N\to\overline{\Q}$ is a non-zero and non-degenerate linear recurrence sequence in $\overline{\Q}$ such that $c:=\max|\al|>1$, where the maximum is taken over all roots $\al$ of $f$. Then for every $\ep>0$ there is an $n_0\in\N$ such that
$$
\forall\,n\ge n_0:\ |f(n)|>c^{(1-\ep)n}\;.
$$
\end{thm}

\noindent
In the full version we will discuss this result in more details in Section 3. Here we only mention that all known proofs are unfortunately non-effective, 
no algorithm is known that would compute for given $f$ and $\ep$ the threshold $n_0$. To state Theorem~\ref{lower_b} in an equivalent 
form we need an important definition.

A {\em power sum} is a formal expression
$$
S(x)=\sum_{i=1}^k p_i(x)\al_i^x
$$
where $x$ is a formal variable, $k\in\N_0$, $p_i\in\overline{\Q}[x]$ are nonzero polynomials, and $\al_i\in\overline{\Q}$ 
are non-zero and pairwise distinct {\em roots of $S(x)$}. $S(x)$ is 
{\em non-empty} if $k\ge1$, and $S(x)$ is {\em non-degenerate} if no ratio $\al_i/\al_j$, $i\ne j$, of two roots is a root of $1$. 
Substituting $x=n\in\N$ makes $S(x)$ to a function $S\cc\N\to\overline{\Q}$. For the {\em empty power sum} with $k=0$ 
this function is the zero sequence. The equivalent form of Theorem~\ref{lower_b} is as follows. 

\begin{thm}
\label{lower_b_ps}
Suppose that $S(x)$ is a non-empty and non-degenerate power sum such that $c:=\max|\al|>1$, where the maximum is taken over all roots $\al$ of $S(x)$. Then for every 
$\ep>0$ there is an $n_0\in\N$ such that
$$
\forall\,n\ge n_0:\ |S(n)|>c^{(1-\ep)n}\;.
$$
\end{thm}

We state the more precise form of Theorem~\ref{main} whose proof we outline in Section 2. Let the tuples in
$${\textstyle
L=\{\overline{a}=(a_0,\ds,a_{2k-1})\in\bigcup_{k\ge0}\Z^{2k}\;|\;a_0\ne0\}
}
$$ 
encode linear recurrence sequences: for $\overline{a}\in L$ let $f(n)=f(\overline{a},n)\cc\N\to\Z$ be the linear recurrence sequence given by
$f(i)=a_{k-1+i}$ for $1\le i\le k$, and by $f(n+k)=\sum_{i=0}^{k-1}a_if(n+i)$ for $n\in\N$. 

\begin{thm}
\label{main_gene}
In Section 2 we describe an algorithm ${\cal A}\cc L\times\N\to\Z$ such that for every $\overline{a}\in L$, the specialization
$$
{\cal A}(\overline{a},n)\cc\N\to\Z
$$
is a PIO formula for the linear recurrence sequence $f(\overline{a},n)$.
\end{thm}

Theorem~\ref{lower_b}, equivalently Theorem~\ref{lower_b_ps}, is a crucial component of the proof of Theorem~\ref{main_gene}. 
As we will discuss in more details in Section 3 of the full version, in 1977 this was still a conjecture. A proof of it 
was given by A.~van der Poorten in a preprint \cite{poor1} in 1981. Earlier no algorithm evaluating linear recurrence sequences and provably 
effective by Definition~\ref{our_defi} was possible, at least not along the lines of Section 2. 

\section{Outline of the proof of Theorem~\ref{main_gene}}

Theorem~\ref{main_gene} follows from the next proposition. Recall that $L$ is the set of integral tuples with even lengths and nonzero first coordinate
encoding linear recurrence sequences. 

\begin{prop}
\label{algor_B}
We describe an algorithm ${\cal B}$ with inputs $\overline{a}$ in $L$ and outputs ${\cal B}(\overline{a})=X$, where $\overline{a}=(a_0,\ds,a_{2k-1})$, $X\sus[m]=\{1,2,\ds,m\}$, $m=m(k)=k^2!$, and
such that the following holds.
\begin{enumerate}
    \item If $j\in X$ then there is a polynomial $q\in\Q[x]$ that is identically zero or has $\deg q<k$ and such that $f(\overline{a},n)=q(n)$ for every $n\equiv j$ modulo $m$.
    \item If $j\in[m]\setminus X$ then there is a $c>1$ and an $n_0\in\N$ such that $|f(\overline{a},n)|>c^n$ for every $n\equiv j$ modulo $m$ with $n\ge n_0$.
\end{enumerate}
\end{prop}

\noindent
{\bf Proof of Theorem~\ref{main_gene}. }The algorithm ${\cal A}$ works as follows. For given $\overline{a}$ in $L$ and $n$ in $\N$ we compute ${\cal B}(\overline{a})=X$ 
in $O(1)$ time. Then we determine, in $\mathrm{poly}(\log(1+n))$ time, if $n$ reduced mod $m=m(k)$ is in $X$ or not. In the former case we find in $O(1)$ time by 
Lagrange's interpolation the polynomial $q(x)$ and compute $f(\overline{a},n)=q(n)$ in $\mathrm{poly}(\log(1+n))$ time by direct evaluation of $q(x)$ at $x=n$. In the latter case 
we compute $f(\overline{a},n)$ in $\mathrm{poly}(n)$ time by the defining recurrence by the simple algorithm. The implicit constants in $O(\cdot)$s and $\mathrm{poly}(\cdot)$s
depend only on $\overline{a}$. It follows, especially from the lower bound in part 2, that this is a PIO formula for $f(\overline{a},n)$. As we remarked above, presently the implicit constant $c$ in the complexity bound for ${\cal A}$ is non-effective. 
\eproof

Proposition~\ref{algor_B} follows from Theorem~\ref{lower_b_ps} and the next seven propositions. Most of them, if not all, 
are well known results from algebra and from the theory of linear recurrence sequences. For the sake of completeness and reader's convenience, 
in the full version we prove all of them, except for Proposition~\ref{kronecker} (which is a standard result); here all proofs are omitted. 
Even in the full version we do not prove two steps in the proof of Theorem~\ref{main_gene}: we omit the proofs of the Subspace Theorem
and of the SML theorem which are both needed to deduce Theorem~\ref{lower_b_ps} in Section 3 of the full version. 

A power series $\sum_{n\ge0}a_n x^n$ with coefficients $a_n\in\Z$
is {\em primitive} if these coefficients are altogether coprime, which means that if $d\in\N$ divides $a_n$ for every $n\in\N_0$ then $d=1$. The next proposition belongs to the circle of results known in
algebra as Gauss' lemma (S.~Lang \cite[Chapter IV.2]{lang}).

\begin{prop}
\label{gauss}
The product of two primitive power series in $\Z[[x]]$ is a primitive power series.
\end{prop}

Let $f\cc\N\to\Q$ be a linear recurrence sequence in $\Q$ with recurrence order $k\in\N_0$. Note that $k$ is minimum for the given $f$ if and only if we have a formal identity in $\Q[[x]]$
$$
\sum_{n\ge1}f(n)x^n=\frac{p(x)}{1-a_{k-1}x-\ds-a_0x^k}
$$
where $p\in\Q[x]$, $p$ is identically zero or has $\deg p\le k$, and $p(x)$ is coprime to the denominator $1-a_{k-1}x-\ds-a_0x^k$ in $\Q[x]$.
The next result is known as Fatou's lemma, after P.~Fatou's article \cite{fato} (another and more widely known Fatou's lemma concerns integration). 
We take the proof from R.\,P.~Stanley's book \cite[p. 275]{stan} where it is attributed to A.~Hurwitz.

\begin{prop}
\label{ZaQ}
If $f\cc\N\to\Z$ is a linear recurrence sequence in $\Q$ then $f$ is a linear recurrence sequence. In fact, the shortest recurrence  for $f$ in $\Q$ is an integral 
recurrence. 
\end{prop}

For any sequence $f\cc\N\to X$ and $m\in\N$ we call the $m$ sequences $f_{m,j}\cc\N\to X$, 
$$
f_{m,j}(n)=f(j+m(n-1))\cc\N\to X\;\text{ with}\;j\in[m]\;,
$$
{\em $m$-sections of $f$}.
\begin{prop}
\label{mjsec}
Let $f\cc\N\to\Z$ be a linear recurrence sequence and $m\in\N$. Then every $m$-section of $f$ is a linear recurrence sequence.
\end{prop}

\noindent
It folows from the proof that all $m$ sequences $f_{m,j}$ satisfy the same recurrence with an order that is at most the order of the recurrence for $f$.

\begin{prop}
\label{algor_C}
We describe an algorithm ${\cal C}$ that for each input 
$$
\overline{a}=(a_0,\,a_1,\,\ds,\,a_{2k-1})\in L
$$ 
outputs a power sum ${\cal C}(\overline{a})=S(x)$ such that $f(\overline{a},n)=S(n)$ holds for every $n\in\N$. Moreover, every root of $S(x)$ is a root of $q_{\underline{a}}(x)$ and every coefficient polynomial $p_i(x)$ in $S(x)$ has degree $\deg p_i<k$.
\end{prop}

\noindent
So every linear recurrence sequence $f\cc\N\to\Z$ is represented as $f(n)=S(n)$ by a power sum $S(x)$. The same is true and the same proof works for linear recurrence sequences in $\Q$, and in $\overline{\Q}$.

The advantage of representing sequences by power sums rather than by linear recurrences is uniqueness, which we show next. One can show that 
if $S(x)$ is a power sum then $S(n)$ is a linear recurrence sequence in $\overline{\Q}$ (hence the equivalence of Theorems~\ref{lower_b_ps} 
and \ref{lower_b}), and that if $S(x)$ is such that $S\cc\N\to\Q$ then $S(n)$ is a linear recurrence sequence in $\Q$, but we will not need these results.

\begin{prop}
\label{ps_jednoz}
If two power sums $S(x)$ and $T(x)$ satisfy $S(n)=T(n)$ for every $n\in\N$ then $S(x)=T(x)$.
\end{prop}

\begin{prop}
\label{ps_root1}
If $S(x)$ is a power sum such that $S\cc\N\to\Z$, if this sequence is a linear recurrence sequence, and if each root $\al$ of $S(x)$ has modulus 
$|\al|\le1$, then each root of $S(x)$ is a root of $1$.
\end{prop}

\begin{prop}[Kronecker's theorem]
\label{kronecker}
If $p\in\Z[x]$ is a monic polynomial such that $p(0)\ne0$ and every root $\al$ of $p(x)$ has modulus $|\al|\le1$, then every root of $p(x)$ is a root of $1$.
\end{prop}

\noindent
{\bf Proof of Proposition~\ref{algor_B}. }The algorithm ${\cal B}$ works as follows. Suppose that a $2k$-tuple $\overline{a}=(a_0,\ds,a_{2k-1})$, 
$a_0\ne0$, of integers is given. We compute by Proposition~\ref{algor_C} a power sum 
$$
S(x)=\sum_{i=1}^l p_i(x)\al_i^x
$$ 
such that $S(n)=f(n)=f(\overline{a},n)$ holds for every $n\in\N$. We set $m=m(k)=k^2!$. For each $j\in[m]$ we compute the power sum
$$
S_{m,j}(x)=\sum_{i=1}^{s_j} q_{i,j}(x)\be_{i,j}^x=\sum_{i=1}^l\al_i^{j-m}p_i(j+m(x-1))(\al_i^m)^x=S(j+m(x-1))\;.
$$
So $S_{m,j}(n)=f_{m,j}(n)$ is the $m$-section of $f(n)$. We compute
$$
X=\{j\in[m]\;|\;s_j=0\mbox{ or for all $i$, }|\be_{i,j}|\le1\}\;.
$$
It follows  from the bounds on the degrees 
$$
\deg_{\Q}(\al_i/\al_j)\le k^2\;\mbox{ and }\; \deg_{\Q}(\al_i)\le k
$$
(by Proposition~\ref{algor_C} each root of $S(x)$ is a root of $q_{\underline{a}}(x)$), from the choice of $m$, and from the 
fact that each root of $S_{m,j}(x)$ is an $m$-th power of some $\al_i$, that each $S_{m,j}(x)$ is non-degenerate and the only root of 
$1$ that may appear among its roots $\be_{i,j}$ is $1$ itself.

Let $j\in X$. Then either $S_{m,j}(x)$ is empty and the claim holds with $q(x)=0$ or by Propositions~\ref{mjsec} and \ref{ps_root1} 
all roots $\be_{i,j}$ are roots of $1$. But the above
observation implies that then $s_j=1$ and $\be_{1,j}=1$. So $S_{m,j}(x)=q_{1,j}(x)$. Since $S_{m,j}(n)=q_{1,j}(n)\in\Z$ for every $n\in\N$, 
we see that $q_{1,j}(x)\in\Q[x]$. 
For any $n\in\N$ with $n\equiv j$ modulo $m$ one has
\begin{eqnarray*}
f(n)&=&f_{m,j}(1+(n-j)/m)=q_{1,j}(1+(n-j)/m)\\
&=&q(n),\ q(x)=q_{1,j}(1+(x-j)/m)\in\Q[x]\;.
\end{eqnarray*}
By Proposition~\ref{algor_C}, $\deg q=\deg q_{1,j}\le\max_i\deg p_i<k$. 

Let $j\in[m]\setminus X$. Then $S_{m,j}(x)$ is non-empty and non-degenerate and the maximum modulus of its roots is larger than $1$. So for
$n\in\N$ with $n\equiv j$ modulo $m$ we have by Theorem~\ref{lower_b_ps} for
$$
f(n)=S_{m,j}(1+(n-j)/m)
$$
the lower bound in part 2.
\eproof

\noindent
This concludes the proof of Theorem~\ref{main_gene}.
\eproof

\noindent
If an effective version of Theorem~\ref{lower_b_ps} were available and we could compute $n_0$ from $\overline{a}$ and $\ep$, algorithms 
${\cal B}$ and ${\cal A}$ could be 
simplified. We could then determine if $j\in X$ just by checking if $f_{m,j}(n)$ agrees with $q(n)$, 
for the $q\in\Q[x]$ computed by Lagrange's interpolation, for sufficiently many $n$. 

For other results on linear recurrence sequences treated from the algorithmic perspective see, for example, 
S.~Almagor, B.~Chapman, M.~Hosseini, J.~Ouaknine and J.~Worrell \cite{alma_al} or J.~Ouaknine and J.~Worrell \cite{posi,ulti_posi}.

\section{The SML (Skolem--Mahler--Lech) theorem and the Subspace Theorem}

The proof of Theorem~\ref{lower_b_ps}, due to P.~Corvaja, uses the following two deep theorems (the latter was a milestone in
number theory and its proof is not easy) which are explained (but, of course, not proved) together with P.~Corvaja's proof 
in the full version of this extended abstract. Only to properly define all terms used in the Subspace Theorem takes some time.

\begin{thm}[the SML theorem, \cite{skol,mahl,lech}]
\label{sml}
For any non-empty and non-degenerate power sum  $S(x)$ the set of zeros
$$
\{n\in\N\;|\;S(n)=0\}
$$
is finite.
\end{thm}

\begin{thm}[the Subspace Theorem, \cite{bomb_gubl}]
\label{sub_thm}
Let $K$ be a number field, $M$ be the set of places on $K$, $|\cdots|_v$ for $v\in M$ be the normalized absolute values, and let $S\sus M$ be a finite subset containing 
all Archimedean places. For every $v\in S$ let 
$$
L_{v,i}(\overline{x})=L_{v,i}(x_1,\ds,x_n),\ i\in[n]\;,
$$
with $n\in\N$ be $n$ linearly independent linear forms with coefficients in $K$, and let $\ep>0$ be given. Then 
all $S$-integral solutions $\overline{x}\in O_{S,K}\setminus\{\overline{0}\}$ of the inequality
$$
\prod_{v\in S}\prod_{i=1}^n|L_{v,i}(\overline{x})|_v<H(\overline{x})^{-\ep}
$$
lie in finitely many hyperplanes in $K^n$. 
\end{thm}

\bigskip\noindent
{\bf Acknowledgements. }I want to thank Michel Waldschmidt for communicating to me P.~Corvaja's derivation of Theorem~\ref{lower_b_ps}
from Theorems~\ref{sml} and \ref{sub_thm}, and to Pietro Corvaja for his kind permission to use his proof.

\end{document}